# INCORPORATION OF GEOMETALLURGICAL MODELLING INTO LONG-TERM PRODUCTION PLANNING

*A. Navarra[1], T. Grammatikopoulos[2], and K. Waters[3]

[1]*Universidad Católica del Norte*
*0610 Angamos*
*Antofagasta, Chile*
*(\*Corresponding author: anavarra@ucn.cl)*

[2]*CSGS Canada Inc.*
*185 Concession Street*
*Lakefield, Ontario, Canada K0L 2H0*

[3]*McGill University*
*3610 University Street*
*Montreal, Quebec, Canada  H3A 0C5*

## ABSTRACT

Strategic decisions to develop a mineral deposit are subject to geological uncertainty, due to the sparsity of drill core samples. The selection of metallurgical equipment is especially critical, since it restricts the processing options that are available to different ore blocks, even as the nature of the deposit is still highly uncertain. Current approaches for long-term mine planning are successful at addressing geological uncertainty, but do not adequately represent alternate modes of operation for the mineral processing plant, nor do they provide sufficient guidance for developing processing options. Nonetheless, recent developments in stochastic optimization and computer data structures have resulted in a framework that can integrate operational modes into strategic mine planning algorithms. A logical next step is to incorporate geometallurgical models that relate mineralogical features to plant performance, as described in this paper.

# INTRODUCTION

Geometallurgy is the analysis of spatially correlated geological data, for the predictive modelling of extractive metallurgical operations (Alruiz et al., 2009; Suazo et. al., 2010; Navarra et al., 2017a). It contributes to a system-wide perspective of the mineral value chain (Navarra et al., 2017b), thus it coordinates mining, stockpiling and blending, with the individual unit operations that occur within mineral processing plants. This perspective overcomes traditional interdisciplinary barriers, leading to more effective strategic decisions and operating practices.

At some mines, the logic employed by managers and engineers may be fundamentally correct, but could benefit from quantitative fine-tuning of operational settings, or by streamlining the incoming data. In particular, geometallurgical models have been developed to relate grinding circuit throughput to incoming mineralogical data (Alruiz et al., 2009), and then to predict flotation kinetics (Suazo et. al., 2010). With recent advances in digitization (Chambers & Thornton, 2016), modern information systems are now a means to integrate geometallurgical models into daily, monthly and longer-term decision-making processes.

Decision-makers must consider a range of data, which is available under different timeframes, and with differing levels of confidence (Lamghari & Dimitrakopoulos, 2016). To mitigate risk, long-term production plans must be sufficiently flexible to allow for optimal short-term decisions, as more detailed information becomes available. Plans should not be overly dedicated to a single possible scenario (e.g. the "mean" scenario). Rather, they should be configured so that they may perform well for the entire distribution of possible scenarios. This is particularly true for mine production, which is subject to various forms of environmental and market uncertainty.

Geological uncertainty is especially crucial for mineral processing plants that are fed by nearby orebodies (Montiel & Dimitrakopoulos, 2015; Goodfellow & Dimitrakopoulos, 2016). This type of uncertainty is related to the sparsity of drill core samples used to characterize the orebody. As depicted in Figure 1, traditional techniques consider only a single geological scenario (i.e. the kriging mean) for long-term production planning, hence they do not measure the confidence of net present value (NPV) estimates, under geological uncertainty. Newer techniques consider a sequence of possible scenarios obtained through conditional simulation (Remy et al., 2009), described in the following section; these techniques employ two-stage stochastic optimization, described in the final section of this paper. In comparison to the traditional deterministic approaches (Figure 1), the current stochastic approaches (Figure 2) have been shown to produce more adaptable mine plans, thereby increasing the expected NPV of mining operations by over 20%, which may correspond to hundreds of millions of dollars (Goodfellow & Dimitrakopoulos, 2016).

An important mechanism to mitigate uncertainty is to develop alternate modes of operation (Navarra et al., 2017b). For instance, a geometallurgical model may determine that a certain type of ore is economical if it undergoes a coarse grind. However, the change from a fine grind to a coarse grind may involve emptying the semi-autogenous grinding mill (SAG mill) and introducing larger steel balls, which would correspond to a change of operational mode. This type of alternation should not be done every day, for instance, but only according to tactical criteria, which considers forecasts from the long-term mine plan, as well as the incoming geometallurgical data. Generally, operational modes can harmonize various segments of the mineral value chain, e.g. an operational mode prescribes milling parameters, while simultaneously prescribing the corresponding upstream parameters for stockpiling and blending, as well as the corresponding downstream parameters for flotation. Thus, operational modes are a fundamental consideration within a system-wide perspective.

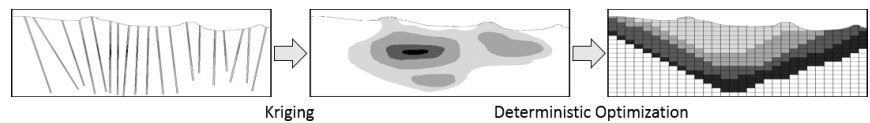

Figure 1. Traditional approach which applies kriging to drill core samples to construct a single-scenario geological model, and deterministic optimization to produce a long-term mine plan

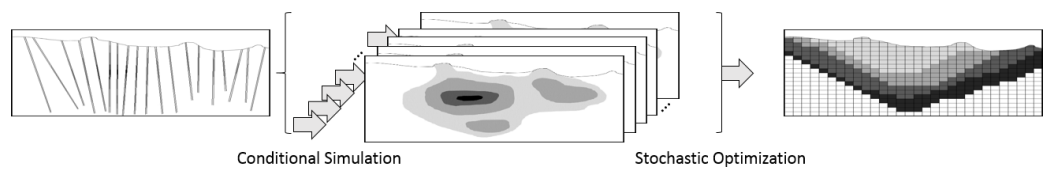

Figure 2. Current approach which applies conditional simulation to construct several geological scenarios, and stochastic optimization to produce a long-term mine plan

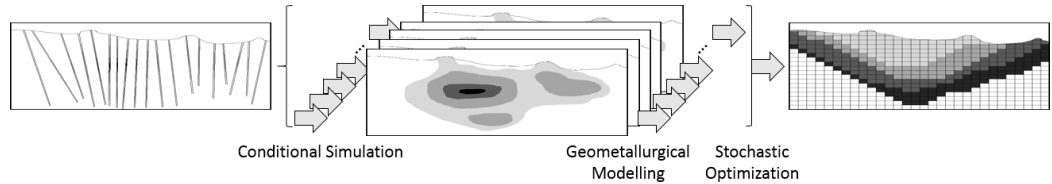

Figure 3. Extension of current approach to long-term mine planning that incorporates geometallurgical modelling

Figure 3 is an extension of Figure 2 that illustrates the incorporation of geometallurgical modelling into long-term mine planning. The additional modelling step is to categorize (or cluster) the blocks in terms of how conducive they are to the proposed operational modes; the resulting categories are known as geometallurgical units (Alruiz et al., 2009; Suazo et. al., 2010). Even in the early stages of mine development, the approach of Figure 3 may be used to determine whether or not an operational mode is economically viable; the stochastic optimization is run with and without the additional mode. Typically, the additional mode would require more processing capacity and/or equipment, and hence a greater capital expenditure. The mode is economically justified if it leads to a statistically significant increase in operational NPV that offsets the additional capital expense (Navarra et al., 2017a).

The establishment of distinct geometallurgical units depends on advanced characterization techniques, such as quantitative evaluation of minerals by scanning electron microscopy (QEMSCAN). As described in the following section, the compositional and morphological features, which can be observed with the QEMSCAN provide insight to determine which combination of blending, milling and concentration parameters might constitute operational modes that are economically viable (Grammatikopoulos et al., 2013; Jordens et al., 2016; Little et al., 2016). Different sources of ore can be regarded as distinct geometallurgical units, depending on whether or not they exhibit categorically different responses to the viable operational modes (Navarra et al., 2017a).

Geometallurgical modelling and stochastic optimization are active areas of research. There is substantial economic incentive to merge these areas in the development of long-term mine plans that mitigate geological uncertainty (Goodfellow & Dimitrakopoulos, 2016). Considering the diverse challenges that are observed at different mine sites, it is important to establish unifying concepts, which serve as a starting point for customized quantitative solutions (Alruiz et al., 2009). Continued collaboration between industry and academia will provide the next generation of tools that employ a system-wide perspective to evaluate and optimize mining projects.

## CONDITIONAL GEOMETALLURGICAL MODELLING

Current mine plan algorithms assume a simple relation between head grade and recovery. For different ore classes, the recovery is often taken to be constant, as is the milling cost per tonne. Indeed, a valid starting point for geometallurgical modelling (transition from Figure 2 to Figure 3) is to determine the milling costs for these classes in order to attain a prescribed recovery; this usually involves empirical hardness-and-throughput studies. Ultimately, the objective of this modelling phase is to determine the economic value and time required to process individual blocks to the proposed set of operational modes.

Geometallurgical models are developed by considering the following four aspects, which are listed in increasing level of detail:

- Mineralogy: Distribution and correlations between mass fractions of different minerals
- Liberation: Distribution of minerals into separate particles as a function of particle size, exposure and degree of freedom
- Texture: Distribution, shape and orientation of crystal grains within the mineral matrix
- Mineral Chemistry: Correlations between the preceding aspects and the presence of elements of interest, or unwanted penalty elements

Certain geological attributes may describe more than one aspect, e.g. a high Phase Specific Interfacial Area (PSIA) may indicate the interlocking of grains (Little et al., 2016), and is descriptive of both liberation and texture; the microstructure of Figure 4b must undergo considerably more grinding than Figure 4a before the valuable (dark) mineral can be liberated and recovered. Together, liberation- and texture-related attributes constitute the mineral morphology ("shape logic").

QEMSCAN has become a standard tool to obtain fundamental data for all four aspects (Grammatikopoulos et al., 2013); it can quantify a diversity of mineralogical attributes: liberation, exposure, grain size, etc. Electron microprobe analysis and laser ablation ICP-MS are used determine the mineral chemistry that is critical with respect to metals of interest. In addition, time of flight secondary ion mass spectrometry (ToF-SIMS) is used to analyze individual atoms or molecules observed on mineral surfaces (Acres et. al., 2010), for example, to determine mineral losses during flotation.

As part of a stochastic framework, compositional and morphological attributes are assigned to each mining block $b$, depending on the geological scenario $g$. Moreover, the value $v_{bg}^{Process}$ of processing block $b$ under scenario $g$ is a function of the operational mode $o$, and may be modelled at a high level as

$$v_{bg}^{Process}(o) = -c_{bg}^{Process}(o) + \sum_{s \in \mathcal{S}} w_{sbg} \cdot r_{sbg}(o) \cdot v_s^{Selling}(o) \qquad (1)$$

$c_{bg}^{Process}(o)$ is the cost per unit mass of processing block $b$ under scenario $g$, $w_{sbg}$ is the weight fraction of species $s$ within the ore of block $b$ under geological scenario $g$, $r_{sbg}(o)$ is the fraction of species $s$ within block $b$ that is recovered into its intended concentrate stream, and $v_s^{Selling}(o)$ is the net selling value of species $s$; $\mathcal{S}$ denotes the set of species under consideration. For the purposes of this paper $v_{bg}^{Process}$, $c_{bg}^{Process}$ and $v_s^{Selling}$ could be denominated in US dollars per metric tonne ($/T).

Long-term mine plans typically span between 5 and 20 years, and must therefore consider the time value of money. This is usually represented through the use of an annual discount rate $d$. For instance, $v_{bgt}^{Process}(o)$ denotes the discounted value of processing block $b$ under scenario $g$ and in time period $t$. If there are $n_T$ periods under consideration, and they each correspond a one-year duration, then

$$v_{bgt}^{Process}(o) = \frac{v_{bg}^{Process}(o)}{(1+d)^{t-1}} \qquad (2)$$

for periods $t = 1, 2, \ldots, n_T$. Equation 2 is such that $v_{bg1}^{Process} = v_{bg}^{Process}$.

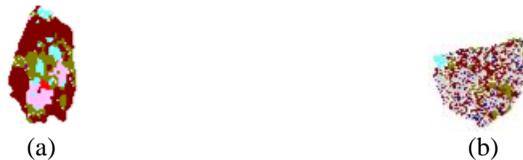

(a)  (b)

Figure 4. Examples of mineral grains having (a) low PSIA and (b) high PSIA.

Complementing the processing value $v_{bg}^{Process}$, is the processing rate $\dot{m}_{bg}^{Process}$, also known as the throughput. $\dot{m}_{bg}^{Process}(o)$ may be taken as the mass rate of block $b$ under geological scenario $g$, if processed under mode $o$, and is measured in T/hour. For simpler models, the processing rate can be fixed by the mode, regardless of the geological attributes, so that $\dot{m}_{bg}^{Process}(o) = \dot{m}^{Process}(o)$. However, more detailed representations of modes may automatically adjust the rate as a function of the geological attributes. For example, a mode may feature dynamic changes in the recirculating load of a grinding circuit, as a function of feed hardness.

Within Equation 1, the processing costs $c_{bg}^{Process}$, recoveries $r_{sbg}$ and selling values $v_s^{Selling}$ are generally taken as functions of the mode $o$, whereas the head grades $w_{sbg}$ are not. More advanced models consider the purity level of concentrate streams as a function of the operational mode $o$, thus affecting the selling values; otherwise, these values could be taken as simple parameters, $v_s^{Selling}(o) = v_s^{Selling}$.

Under a given mode of operation $o$, Equation 1 assumes that species have value only if they are recovered into their intended concentrate stream. This is not applicable, however, when a valuable species is recovered into several concentrate streams, presumably with different selling values. For example, at certain copper-gold mines, some of the gold is recovered by gravity separation, while a large portion passes through the flotation circuit into the copper sulphide concentrate, to be recovered after smelting. In this case, gravity recovered gold (GRG) might have a higher selling value than the flotation recovered gold (FRG), since it requires less downstream processing. Equation 1 could still be applicable if GRG and FRG are somehow regarded as a different species. Alternatively, $r_{Au,bg}(o) \cdot p_{Au}(o)$ could be replaced by a single (more complex) function of the form $f_{Au,bg}(o)$ that would better describe the partitioning of gold between the two streams.

In most situations, valuable species are (primarily) directed toward a single concentrate stream. Thus, the following simplified version of Equation 1 is a practical framework for geometallurgical modelling:

$$v_{bg}^{Process}(o) = -c_{bg}^{Process}(o) + \sum_{s \in \mathcal{S}} w_{sbg} \cdot r_{sbg}(o) \cdot v_s^{Selling} \qquad (3)$$

in which the net selling values $v_s^{Selling}$ are assumed to be independent of the mode. From this, hypothesis testing and multilinear regression are used to develop empirical expressions for $c_{bg}^{Process}(o)$ and $r_{sbg}(o)$, as well as for $\dot{m}_{bg}^{Process}(o)$, that incorporate critical compositional and morphological attributes (Little et al., 2016). Additional fundamental knowledge about the kinetics of individual unit operations may help identify which terms could be modelled as nonlinear, resulting in semi-empirical expressions.

Equations 1 and 3 do not consider the cost or time of changing from one mode of operation to another. Ideally, the changing occurs relatively infrequently, so that these costs and times are negligible compared to the capital expenditures and the life of the mine. Indeed, the most drastic mode changes could be synchronized with planned maintenance shut-downs. Furthermore, stockpiling policies can be formulated to provide medium- and short- term flexibility in the timing of operational changes (Navarra et al., 2017b).

Assuming that there is enough short-term flexibility, long-term mine plans may assume that each ore block is processed only according to the value-throughput-maximizing operational mode that is available for it (Goodfellow & Dimitrakopoulos, 2016), depending on which geological scenario arises. Supposing that a proposed equipment configuration can support a set of operational modes $\mathcal{O}$, the value-throughput-maximizing mode to process $b$ under geological scenario $g$, is denoted

$$o_{bg} = \underset{o \in \mathcal{O}}{\operatorname{argmax}} \left( v_{bg}^{Process}(o) \cdot \dot{m}_{bg}^{Process}(o) \right) \qquad (4)$$

A different selection of metallurgical equipment generally implies different operational modes $\mathcal{O}$, which impacts the long-term mine plan.

In practice, the individual scenarios *g* that are considered in Equations 1–4 are generated according to geostatistical conditional simulation (Remy et al., 2009); these techniques constitute a kind of Monte Carlo simulation, due to their use of random number generation. Beginning with the drill core data, the simulation fills in the unknown blocks using randomly generated data. In principal, each scenario has an equal probability of occurring, although they are notably different from each other, since they are each generated from a different set of random numbers. In this context "conditional" implies that the simulations preserve local spatial averages and local spatial variances, i.e. first and second order conditions, respectively. Typically, 10 to 20 scenarios provide a sufficient representation of geological uncertainty, as each additional scenario has a diminishing impact on the resulting mine plan (Montiel & Dimitrakopoulos, 2015; Goodfellow & Dimitrakopoulos, 2016).

Conditional simulation is an area of ongoing research. However, certain techniques have come into relative provenance, and are commonly available as open-source software (Remy et al., 2009). This is especially true for the two-point techniques (Figure 5), including Sequential Gaussian Simulation (SGS) and Sequential Indicator Simulation (SIS). The former is appropriate when the attribute is parametric (i.e. Gaussian-transformable), whereas the latter for non-parametric attributes, including discretized and categorical attributes. SGS is especially useful for numerical attributes that vary continuously over space, including common mineral grades. However, attributes that have strongly correlated extreme values, or high values that are correlated differently than low values, should be discretized into geological categories, and subject to SIS. This is the case for mineral deposits, in which rich zones are oriented preferentially, and are superimposed over a more homogeneous background.

SGS, SIS and other two-point techniques employ so-called variogram models that describe how the attributes may vary between two sampling points, as a function of their distance. The two-point approach is appropriate when the region is relatively free of organized macroscopic structures (Remy et al., 2009). Otherwise, multiple-point simulations can be formulated to capture increasingly intricate macroscopic patterns, as observed in gold and diamond deposits the most used multi-point method is Single Normal Equation Simulation (SNES). Figure 6 is an illustrative comparison between two-point and multi-point simulations.

Of particular relevance for geometallurgical modelling, is the incorporation of so-called secondary attributes, especially with regard to liberation, texture and mineral chemistry. Sequential Gaussian Co-Simulation (SGCS) and Sequential Indicator Co-Simulation (SGCS) are extensions of SGS and SIS, respectively, which simulate a sparsely sampled attribute (i.e. the primary attribute), benefitting from its correlations with more densely sampled secondary attributes (Remy et al., 2009). (Similar extensions also exist for multi-point techniques). For example, there may be extensive data for a spatial distribution of galena, but sparse data for the mercury uptake within the galena; SGCS conditions sequentially simulate blocks using galena spatial correlation, as well as galena-mercury correlation and mercury correlation, thus providing more reliable estimates than SGS. This use of conditional simulation, in combination with the operational modelling efforts of Equations 1 – 4, may be called "conditional geometallurgical modelling".

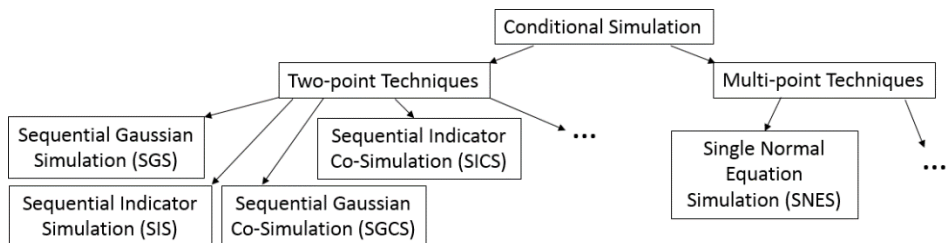

Figure 5. Categorization of selected conditional simulation techniques

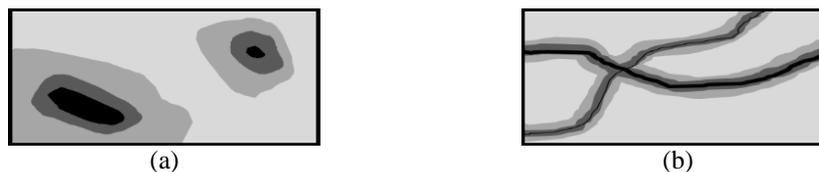

(a)            (b)

Figure 6. Depiction of simulation from (a) two-point, and (b) multi-point simulations

Characterization data from automated mineralogical techniques (e.g. QEMSCAN) lead to the formulation of basic hypotheses to identify which operational modes $o$ may be well-suited for the observed compositional and morphological attributes, and subsequently, which of these attributes will form significant terms within the functions for processing values $v_{bg}^{\text{Process}}(o)$ and rates $\dot{m}_{bg}^{\text{Process}}(o)$. However, some consideration must be given to the conditional simulation techniques that are used to generate the underlying geological scenarios. Ultimately, these functions and processing options form a mapping between the geostatistical scenarios and the geometallurgical units (Figure 3), which are used for long-term mine planning.

## TWO-STAGE MINE PLAN OPTIMIZATION

Two-stage optimization generally considers both a long- and a short-term timeframe. In particular, long-term mine plans are subject to a range of possible scenarios (Figures 2 and 3), based on multi-year forecasts of highly uncertain geological attributes. Subsequently, processing decisions are finalized within a short-term timeframe, benefitting from monthly and weekly forecasts that are more precise.

For a deposit that is described by a set of blocks $\mathcal{B}$, a long-term mine plan that spans $n_T$ time periods may be represented as $x = \{\mathcal{B}_{xt}\}_{t=1}^{n_T}$ in which $\mathcal{B}_{xt} \subset \mathcal{B}$ is the set of blocks that are to be mined in time period $t$. To represent short-term processing decisions, $y_g = \{\{y_{gbt}\}_{b\in\mathcal{B}}\}_{t=1}^{n_T}$ is a matrix of fractional values such that $y_{gbt}$ is the portion of block $b$ that is processed in period $t$. In summary:
- $x$ describes when the blocks are *mined* (long-term decisions)
- $y_g$ describes when the blocks are *processed* (short-terms decision if scenario $g$ is realized)

This formulation of $x$ is consistent with typical mine planning algorithms (Lamghari & Dimitrakopoulos, 2016), supposing that a block $b$ is either entirely mined in period $t$ (i.e. $b \in \mathcal{B}_{xt}$), or not at all (i.e. $b \notin \mathcal{B}_{xt}$), hence a discrete decision. However, the short-term variables $y_{gbt}$ are continuous, allowing the blocks to be processed and/or stockpiled over several periods. Nonetheless,

$$(b \in \mathcal{B}_{xt}) \Rightarrow (y_{gbt'} = 0 \text{ for all preceding periods } t' < t \text{ and all scenarios } g) \quad (5)$$

meaning that a block cannot be processed prior to being mined.

The objective of stochastic long-term mine planning is to maximize the expected NPV,

$$\mathbb{E}[\text{NPV}](x) = -\frac{1}{n_G}\sum_{g=1}^{n_G}\sum_{t=1}^{n_T} c_{gt}^{\text{Mining}}(x) + \frac{1}{n_G}\sum_{g=1}^{n_G}\left(\max_{y_g \in \mathcal{Y}_g(x)} \sum_{t=1}^{n_T} v_{gt}^{\text{Process}}(x, y_g)\right) \quad (6)$$

simultaneously accounting for $n_G$ geological scenarios, in which $c_{gt}^{\text{Mining}}(x)$ is the discounted mining cost in period $t$ under scenario $g$ as incurred by long-term mine plan $x$, and $v_{gt}^{\text{Process}}(x, y_g)$ is the discounted value obtained in period $t$ under scenario $g$, as incurred by long-term plan $x$ and short-term decisions $y_g$. Due to Equation 5, the set of feasible short-term decisions $\mathcal{Y}_g(x)$ is indeed a function of the mine plan $x$; short-term feasibility will be discussed further below.

Within Equation 6, the first double-summation (expected mining costs) and the second double-summation (expected processing values) are both taken as averages over the $n_G$ geological scenarios. Moreover, the second double summation has an embedded maximization operation; the outer optimization of $x$ with respect to [NPV] thus incorporates an inner optimization stage that determines the short-term processing decisions $y_g$, as more precise geometallurgical information becomes available.

PREPRINT – (SUBMITTED TO THE INTERNATIONAL JOURNAL OF MINERAL PROCESSING)

Regarding the first double-summation in Equation 6, mining costs are generally assumed to have the following form (Lamghari & Dimitrakopoulos, 2016):

$$c_{gt}^{\text{Mining}}(x) = \sum_{b \in \mathcal{B}_{xt}} c_{bgt}^{\text{Mining}} \cdot m_{bg} \quad (7)$$

in which $c_{bgt}^{\text{Mining}}$ is the discounted mining cost per unit mass of block $b$ under scenario $g$ in period $t$, and $m_{bg}$ is the total mass of block $b$ under scenario $g$. Thus, every block $b \in \mathcal{B}_{xt}$ contributes a single term $c_{bgt}^{\text{Mining}} m_{bg}$ to the mining costs projected for time $t$, under each scenario $g$. This form permits a reordering of the first double summation in Equation 6, so that

$$\frac{1}{n_G} \sum_{g=1}^{n_G} \sum_{t=1}^{n_T} c_{gt}^{\text{Mining}}(x) = \frac{1}{n_G} \sum_{g=1}^{n_G} \sum_{t=1}^{n_T} \sum_{b \in \mathcal{B}_{xt}} c_{bgt}^{\text{Mining}} \cdot m_{bg}$$

$$= \frac{1}{n_G} \sum_{t=1}^{n_T} \sum_{b \in \mathcal{B}_{xt}} \sum_{g=1}^{n_G} c_{bgt}^{\text{Mining}} \cdot m_{bg}$$

$$= \sum_{t=1}^{n_T} \sum_{b \in \mathcal{B}_{xt}} \left( \frac{1}{n_G} \sum_{g=1}^{n_G} (c_{bgt}^{\text{Mining}} \cdot m_{bg}) \right) = \sum_{t=1}^{n_T} c_t^{\text{Mining}}(x)$$

in which

$$c_t^{\text{Mining}}(x) = \sum_{b \in \mathcal{B}_{xt}} \left( \frac{1}{n_G} \sum_{g=1}^{n_G} c_{bgt}^{\text{Mining}} \cdot m_{bg} \right) \quad (8)$$

is the expected discounted mining cost in period $t$ incurred by plan $x$.

The objective function (Equation 6) can therefore be restated as

$$\mathbb{E}[\text{NPV}](x) = -\sum_{t=1}^{n_T} c_t^{\text{Mining}}(x) + \frac{1}{n_G} \sum_{g=1}^{n_G} \left( \max_{y_g \in \mathcal{Y}_g(x)} \sum_{t=1}^{n_T} v_{gt}^{\text{Process}}(x, y_g) \right) \quad (9)$$

Thus the expected value operator ($\frac{1}{n_G} \sum_{g=1}^{n_G}$) has been "folded into" the $c_{gt}^{\text{Mining}}$ function, thereby cancelling out the $g$ subscript. Unfortunately, a similar approach cannot be applied to eliminate the $g$ from $v_{gt}^{\text{Process}}$, because the expected value operator ($\frac{1}{n_G} \sum_{g=1}^{n_G}$) does not generally commute with the maximization operator ($\max_{y_g \in \mathcal{Y}_g(x)}$), i.e. the average of the maximum does not equal the maximum of the average.

The computationally relevant form of Equation 9 is given by

$$\mathbb{E}[\text{NPV}](x) = -c^{\text{Mining}}(x) + \sum_{g=1}^{n_G} v'^{\text{Process}}_g(x, y_g(x)) \quad (10)$$

in which a single computational variable is used to store the total discounted mining costs, $c^{\text{Mining}}(x) = \sum_{t=1}^{n_T} c_t^{\text{Mining}}(x)$. Also, $n_G$ computational variables are used to store each of the scenario-wise contributions to the discounted processing value, $v'^{\text{Process}}_g(x, y_g(x)) = \frac{1}{n_G} \sum_{t=1}^{n_T} v_{gt}^{\text{Process}}(x, y_g(x))$, in which

$$y_g(x) = \operatorname*{argmax}_{y_g \in \mathcal{Y}_g(x)} v'^{\text{Process}}_g(x, y_g) \quad (11)$$

represents the optimal short-term decisions, subject to scenario $g$ and long-term plan $x$.

Mine planning algorithms consider a set a feasible mine plans $\mathcal{X}$. Herein, a mine plan $x$ is feasible (i.e. $x \in \mathcal{X}$), if it respects the mine tonnage capacity $m_t^{\text{MiningCap}}$ for each period $t$, as well as structural mechanic constraints. In the case of open-pit mining, the structural mechanic constraints include a maximum permissible pit angle (Lamghari & Dimitrakopoulos, 2016); more elaborate constraints are considered for underground mining (Carpentier et al., 2016).

Some caution must be taken when defining the mine tonnage capacity within a stochastic context. A common approach is to treat $m_t^{\text{MiningCap}}$ as an upper-bound on the expected mining tonnage (Lamghari & Dimitrakopoulos, 2016), thereby ignoring the variation in the block masses, i.e. taking $m_{bg} = m_b$ to be independent of the scenario. Another approach is implement an upper bound $p_t^{\text{MiningCap}}$ on the probability of violating the mining capacity,

$$\sum_{g=1}^{n_G} \mathbb{I}\left[\left(\sum_{b \in \mathcal{B}_{xt}} m_{bg}\right) > m_t^{\text{MiningCap}}\right] \le p_t^{\text{MiningCap}} n_G \qquad (12)$$

in which the indicator operator $\mathbb{I}[.]$ is 1 if the statement inside the square brackets is true, and is 0 otherwise; thus the left side of Equation 12 counts the number of scenarios for which the mining capacity is violated, to verify that no more than $p_t^{\text{MiningCap}} n_G$ of the $n_G$ scenarios are in violation.

Starting with an initial solution $x$, mine planning algorithms, including the Bienstock-Zukerberg algorithm, and various metaheuristic approaches such as Simulated Annealing, Tabu Search, Variable Neighbourhood Descent, and Diversified Local Search (Lamghari & Dimitrakopoulos, 2016), execute numerous iterations to find alterations in $x$ that may increase $\mathbb{E}[\text{NPV}]$, while respecting the constraints of $\mathcal{X}$. These algorithms differ with respect to the sequence of alterations that they test, and type of constraints that they can support. In general, the iterations involve transferring blocks into and out of the mining periods, while efficiently updating the components of Equation 10, including $c^{\text{Mining}}(x)$ and $v_g'^{\text{Process}}(x, y_g^*(x))$.

The updating of mining costs $c^{\text{Mining}}(x)$ is direct. For every block $b$ that is introduced into $t$,

$$(b \text{ is introduced into } \mathcal{B}_{xt}) \Rightarrow \left(c^{\text{Mining}} \text{ is increased by } \left(\frac{1}{n_G}\sum_{g=1}^{n_G} c_{bgt}^{\text{Mining}} \cdot m_{bg}\right)\right) \qquad (13)$$

which agrees with Equation 8. Similarly, for every block $b$ that is removed from period $t$,

$$(b \text{ is removed from } \mathcal{B}_{xt}) \Rightarrow \left(c^{\text{Mining}} \text{ is decreased by } \left(\frac{1}{n_G}\sum_{g=1}^{n_G} c_{bgt}^{\text{Mining}} \cdot m_{bg}\right)\right) \qquad (14)$$

The rapid updating of $c^{\text{Mining}}(x)$ requires that the quantities $\left(\frac{1}{n_G}\sum_{g=1}^{n_G} c_{bgt}^{\text{Mining}} m_{bg}\right)$ be computed and saved *a priori*, prior to the first iteration, for each block $b \in \mathcal{B}$, scenario $g \in \{1,2,..,n_G\}$ and period $t \in \{1,2,…,n_T\}$.

The updating of $v_g'^{\text{Process}}(x, y_g(x))$ is more complex. At each iteration, $y_g$ must be repeatedly reoptimized while respecting the short-term constraints of $\mathcal{Y}_g(x)$. In addition to the mining-processing precedence (Equation 5), $\mathcal{Y}_g(x)$ enforces a maximum processing capacity in all periods $t$,

$$\sum_{t'=1}^{t}\sum_{b \in \mathcal{B}_{xt'}} \left(\frac{m_{bg}}{\dot{m}_{bg}^{\text{Process}}}\right) y_{gbt} \le d_t^{\text{ProcessCap}} \qquad (15)$$

in which $d_t^{\text{ProcessCap}}$ is the duration that the processing systems are available during period $t$. In a simplified case, there may be a single processing rate $\dot{m}^{\text{Process}}$ that is applied to all blocks and all scenarios so that $\dot{m}_{bg}^{\text{Process}} = \dot{m}^{\text{Process}}$; Equation 15 could then expressed in terms of a mass capacity $m_t^{\text{ProcessCap}} = \dot{m}^{\text{Process}} d_t^{\text{ProcessCap}}$, as is presumed by current techniques (Lamghari & Dimitrakopoulos, 2016). Otherwise,

$$\mathcal{B}_{xgt} = \{b_{(xgt)1}, b_{(xgt)2}, \cdots, b_{(xgt)(i-1)}, \overset{\text{Cutoff block}}{\underset{\downarrow}{b_{(xgt)i}}}, b_{(xgt)(i+1)}, \cdots, b_{(xgt)n_{xt}}\}$$

Figure 7. Pointer which identifies the cutoff block within the ordering $\mathcal{B}_{xgt}$, and can be shifted left or right as the result of an optimization iteration

the time-based capacity $d_t^{\text{ProcessCap}}$ provides a common basis to tally different blocks that are to be processed at different rates.

The efficient updating of $v'^{\text{Process}}_g(x, y_g(x))$ under the restrictions of Equations 5 and 15 was recently accomplished using customized data structures (Navarra et al., 2016; Navarra et al., 2017a). Previous approaches were only able to approximately satisfy the plant capacity constraint (Equation 15), and could not support realistic representations of alternate modes of operation, and were therefore ill-suited to support geometallurgical models (Figure 3).

For each block set $\mathcal{B}_{xt}$, the new approach maintains $n_G$ orderings, one for each geological scenario $g$. Each of these orderings contains all of the blocks of $\mathcal{B}_{xt}$, and is of the form $\mathcal{B}_{xgt} = \{b_{(xgt)1}, b_{(xgt)2}, \ldots, b_{(xgt)\text{nxt}}\}$; these orderings each constitute a priority list to decide which blocks should undergo immediate processing versus being stockpiled. The orderings are therefore sorted from most to least value throughput,

$$v^{\text{Process}}_{b_{(xgt)ig}} \dot{m}^{\text{Process}}_{b_{(xgt)ig}} \geq v^{\text{Process}}_{b_{(xgt)jg}} \dot{m}^{\text{Process}}_{b_{(xgt)jg}} \quad (16)$$

for all $b_{(xgt)i}, b_{(xgt)j} \in \mathcal{B}_{xgt}$ such that the ranking index $i$ is smaller than ranking index $j$, (i.e. $i < j$). Indeed, each scenario provides a different ordering (priority list) for the processing. Moreover, these orderings are implemented as red-black trees (Navarra et al., 2016), hence the $n_G$ orderings are efficiently maintained at each iteration, even as blocks are being transferred between periods.

The orderings $\mathcal{B}_{xgt}$ provide an efficient mechanism to update $y_g(x)$, as well as $v'^{\text{Process}}_g(x, y_g(x))$. As illustrated in Figure 7, a pointer is used to identify the cutoff block. Blocks that are to the left of the pointer (i.e. blocks having a high value throughput) are immediately and entirely processed $y_{gbt} = 1$, whereas blocks to the right of the pointer (i.e. blocks having a low value throughput) are completely differed $y_{gbt} = 0$, possibly to be stockpiled and processed in a later periods; the cutoff block may be partially or entirely processed, $0 < y_{gbt} \leq 1$, so as to occupy the remaining processing capacity $d_t^{\text{ProcessCap}}$. As blocks are transferred between processing periods, the pointers are shifted left to restore short-term feasibility (Equation 15), or right to restore short-term optimality (Equation 11), effectively updating $y_g(x)$. Furthermore, the sorted nature of $\mathcal{B}_{xgt}$ ensures that the correct processing values $v^{\text{Process}}_{bgt}$ are removed/added from $v'^{\text{Process}}_g(x, y_g(x))$ as the pointer is shifted leftward/rightward.

The new implementation considers that $v^{\text{Process}}_{bgt}$ and $\dot{m}^{\text{Process}}_{bg}$ are evaluated using the value-throughput-maximizing mode $o_{bg}$, described by Equation 4. Therefore different process configurations, hence different sets of operational modes $\mathcal{O}$, produce different NPV results. Indeed, the data structures described by Equation 16 and Figure 7 have resulted in a computationally efficient framework to incorporate geometallurgical modelling into long-term stochastic mine planning; sample computations have been prepared by Navarra et al. (2017a). Nonetheless, this implementation does not give rigorous consideration of modes that are faster but less profitable than $o_{bg}$, and which may be preferable to stockpiling; this would be an additional level of complexity that requires further advances in algorithms and data structures.

## CONCLUSIONS AND FUTURE WORK

Two-stage optimization algorithms can be used to evaluate alternate operational modes, hence different plant configurations, under conditions of geological uncertainty. Nonetheless, the implementation of these algorithms requires specialized data structures, and is an ongoing area of research. Innovative data structures allow more detailed representation of the mineral processing operations, while maintaining the computational efficiency of the algorithms. The current framework is now able to support geometallurgical models, as they relate to different modes of operation, and processing capacities.

Future development should address the tradeoff between stockpiling a block of material, and processing it through secondary modes of operation (the current framework only considers the value-throughput-maximizing mode for each block). In addition to geological scenarios, future work can also consider several economic scenarios, thereby addressing price uncertainty, as well as geological uncertainty. Other frameworks have been able to implement blending constraints (Goodgellow & Dimitrakopoulos, 2016), which is be compatible with the current framework, and is are especially relevant for the management of penalty elements. It should also be mentioned that stochastic open-pit mine planning algorithms are considerably more advanced than underground planning algorithms (Carpentier et al., 2016); the latter is an especially important area of future research.